\documentclass[12pt]{amsart}
\usepackage{amsmath,amsthm,amsfonts,amssymb,eucal}

\renewcommand {\a}{ \alpha }

\newcommand{\vare}{\varepsilon}

\newcommand{\G}{\Gamma}

\renewcommand{\d}{\delta}

\renewcommand{\t}{\theta}

\newcommand{\T}{\Theta}
\newcommand{\p}{\partial}

\newcommand{\Om}{\Omega}

\newcommand{\R}{ \mathbb R}
\newcommand{\C}{ \mathbb C}
\newcommand{\N}{ \mathbb N}

\newcommand {\GM}{\mathfrak M}

\newcommand {\BC}{\mathbf C}

\newcommand {\BG}{\mathbf G}

\newcommand {\BT}{\mathbf T}

\newcommand{\lu}{\langle}
\newcommand{\ru}{\rangle}

\newcommand{\CV}{\mathcal V}

\newcommand{\CP}{\mathcal P}

\newcommand{\CC}{\mathcal C}

\newcommand{\CE}{\mathcal E}

\newcommand{\tG}{\textsf G}

\newcommand{\tL}{\textsf L}

\newcommand{\tP}{\textsf P}

\newcommand{\tS}{\textsf S}

\newcommand{\tW}{\textsf W}

\newcommand {\bpr}{\boldsymbol\prime}

\DeclareMathOperator{\res}{\restriction}

\newtheorem{thm}{Theorem}[section]

\newtheorem{lem}[thm]{Lemma}
\newtheorem{prop}[thm]{Proposition}

\theoremstyle{definition}

\theoremstyle{remark}

\numberwithin{equation}{section}

\newcommand{\thmref}[1]{Theorem~\ref{#1}}

 \DeclareMathOperator {\Step}{Step}
\DeclareMathOperator {\rank}{rank} 
\DeclareMathOperator {\comb}{comb} \DeclareMathOperator {\Int}{Int}
\newcommand {\wt}{\widetilde} \newcommand {\wh}{\widehat}

\begin{document}

\title[Approximation on metric graphs] 
{On approximation of functions from Sobolev spaces on metric graphs}
\author[M. Solomyak] 
{Michael Solomyak} 
\address 
{Department of Mathematics 
\\ Weizmann Institute\\ Rehovot\\ Israel}
\email{solom@wisdom.weizmann.ac.il} 
\thanks
{Partially supported by the NATO grant PST.CLG.978694.}
\date{2 May 2002}
\subjclass[2000] {Primary 46E39; Secondary 41A45, 41A46.} 
\dedicatory
{Dedicated to Y. Brudny, outstanding mathematician and an old friend} 
 
\begin{abstract} 
Some results on the approximation of functions from the Sobolev spaces
on metric graphs by step functions are obtained. In particular, we show 
that the approximation numbers $a_n$ of the embedding operator of the Sobolev
space $\tL^{1,p}(\BG)$ on a graph $\BG$ of finite length $|\BG|$ into
the space $\tL^p(\BG,\mu)$, where $\mu$ is an arbitrary finite Borel measure
on $\BG$, satisfy the inequality 
\begin{equation*}
a_n\le |\BG|^{1/{p'}}\mu(\BG)^{1/p}n^{-1},\qquad 1<p<\infty.
\end{equation*}
The estimate is sharp for any $n\in\N$.

\end{abstract} 
\maketitle 

\section{Introduction}

A metric graph is a graph whose edges are viewed as non-degenerate
line segments, rather than pairs of vertices as in the case of
the standard (combinatorial) graphs. This difference is reflected
in the nature of functions on the corresponding graph. For a combinatorial
graph this is just a family of numbers $\{f(v)\}$ where the argument $v$ 
runs over the set of all vertices, while a function on
a metric graph is a family of functions on its edges, usually subject
to some matching conditions at the vertices. 

Sobolev spaces $\tL^{1,p}$ on a metric graph $\BG$ are 
defined in a natural way, 
by analogy with their counterparts for a single interval. The local 
properties of
functions from these spaces outside the vertices are evidently the same
as for the case of an interval. However, the global properties may depend
on the geometry of a given graph. We
establish some results on approximation of functions
from $\tL^{1,p}$ by step functions.
The estimates obtained are uniform with respect to all graphs of a fixed 
length and do not depend on the structure of the graph. The 
estimates are sharp with respect to all the parameters involved.
We believe that such results are useful for better understanding of function
spaces on graphs.
\bigskip

An important phase in the development of analysis on metric graphs
was started by the paper \cite{EH} by W. D. Evans and D. D. Harris. Embeddings
of the Sobolev spaces $\tW^{1,p}(\Om)$ in $\tL^p(\Om)$ were studied there
for a wide class of 
domains with irregular boundary. A characteristic
feature of these domains is that they have a ``ridge'', 
this being a metric tree. In \cite{EH} the study of such embeddings was
reduced to the 
investigation of the behavior of the approximation numbers for the weighted
Hardy-type integral operators on the ridge. For $p=2$ approximation numbers
coincide with the singular numbers, and the problem
can be reformulated in terms of the eigenvalue behavior for the ``weighted
Laplacian'' on the tree. From this point of view
the question was analyzed in \cite{NS}. Eigenvalue estimates for the weighted
Laplacian were
obtained there in terms of appropriate partitions of the given tree into
a family of segments. Some of the results of \cite{NS} were considerably
refined by W. D. Evans, D. D. Harris, and
J. Lang in their recent paper \cite{EHL}. The main novelty of \cite{EHL} 
consisted in replacement
of segments, as elements of a partition, with arbitrary compact subtrees.
A thorough analysis 
of the partitions appearing in the process of
approximation allowed the authors to obtain important results for
arbitrary $p,\ 1\le p\le\infty$. In particular, for $p\in(1,\infty)$
they established a Weyl-type asymptotic formula for 
the approximation numbers.
\vskip0.2cm
Our goal in this paper is to consider arbitrary graphs, rather than 
only
the trees. The language of Hardy-type integral operators is no more relevant,
since such operators are well defined only on trees. Instead, we study
embeddings of Sobolev spaces on the graph
$\BG$ into the space $\tL^\infty(\BG)$
and into the spaces $\tL^p(\BG,\mu)$ where $\mu$ is an arbitrary
Borel measure on $\BG$. The character of the results 
obtained makes it apparent that this language is adequate. Following the idea
of \cite{EHL}, we use partitions of a given graph into subgraphs, however
the way of this usage differs from the one in \cite{EHL}. We restrict ourselves
to the case of compact graphs, since the passage to non-compact ones
can be carried out exactly as in \cite{EHL} and does not require new
ideas, as soon as one is interested only in the estimates but not in 
asymptotics.
\vskip0.2cm
Introduce some necessary notations.
Let $\BG$ be a connected graph with the set of vertices 
$\CV=\CV(\BG)$ 
and the set of edges $\CE=\CE(\BG)$. Compactness of a graph
means that $\#\,\CE<\infty$ and hence, also
$\#\,\CV<\infty$. The distance
$\rho(x,y)=\rho_\BG(x,y)$ between any two points $x,y\in\BG$ 
(and thus, the metric
topology on $\BG$), and also the measure $dx$ on $\BG$ 
are introduced in a natural way; see Section 2 for detail. Below
$|E|=|E|_\BG$ stands for the measure of a measurable set $E\subset\BG$. If in 
particular $E=e$ is an edge, then $|e|$ is its length. 

\bigskip

Below the symbol $\GM(\BG)$ stands for the set of all finite Borel measures on
$\BG$. For $1\le p\le\infty$, we denote by $\|\cdot\|_{p,\mu}$
the norm in the space $\tL^p(\BG,\mu)$, i.e.
\begin{equation*}
\|u\|_{p,\mu}=\|u\|_{\tL^p(\BG,\mu)}=
\biggl(\int_{\BG} |u|^pd\mu\biggr)^{1/p},\qquad p<\infty,
\end{equation*}
with the standard change if $p=\infty$. If the measure $\mu$ is
absolutely continuous, i.e. $d\mu=Vdx$, then we write $V$ instead of $\mu$ 
in the above notations. We drop the index $\mu$ (or $V$) if $d\mu=dx$.

A function $u$ on $\BG$ belongs to
the Sobolev space $\tL^{1,p}=\tL^{1,p}(\BG)$, 
 if $u$ is 
continuous on $\BG$ and its restriction to each edge $e$ 
has the distributional derivative $u'$ which is a function from $\tL^p(e)$. 
The functional $\|u'\|_{\tL^p(\BG)}$
defines on $\tL^{1,p}$ a semi-norm vanishing on the one-dimensional
subspace of constant functions.

We say that $v$ is a {\sl step function} on $\BG$ and write $v\in\Step(\BG)$,
if $v$ takes only a finite
number of different values, each one on a connected subset of $\BG$. 
Any function $v\in\Step(\BG)$ can be represented as a
linear combination of characteristic functions of mutually disjoint
connected subsets. We write $v\in\Step_n(\BG)$, if for $v$ there exists
a representation with the number of terms less or equal to $n$.

We are interested in the approximation of functions 
$u\in\tL^{1,p} (\BG)$
by functions $v\in\Step_n(\BG)$. More exactly, we study two problems:
the uniform
approximation (i.e., approximation in the metric $\|\cdot\|_\infty$)
and approximation in the metric $\|\cdot\|_{p,\mu}$. 
In the first problem we construct a mapping 
$Z_p:\tL^{1,p} (\BG)\to\Step_n(\BG)$ such that 
$\|u-Z_p u\|_\infty\le C_p(\BG)(n+1)^{-1}\|u'\|_p$. This problem is elementary 
for $p=\infty$, when the operator $Z_\infty$ can be chosen linear and
$C_\infty(\BG)=|\BG|$.
For $p<\infty$ a linear mapping $Z_p$ with the required properties
does not exist but we find a 
non-linear mapping
which gives the same rate of approximation, with
$C_p(\BG)= |\BG|^{1/{p'}}$. In the second
problem we establish a similar result by means of a linear approximation
operator; this operator depends on the measure $\mu$.

\bigskip

Below we present formulations of the typical results. 
\begin{thm}\label{1:uni}
Let $\BG$ be a compact graph and $1\le p\le\infty$. Then for any function
$u\in\tL^{1,p} (\BG)$ and any $n\in\N$ there exists a function
$v\in\Step_n(\BG)$ such that
\begin{equation}\label{1:1}
\|u-v\|_\infty\le\frac{|\BG|^{1/{p'}}\|u'\|_p}{n+1}.
\end{equation}
If $p=\infty$, the mapping $u\mapsto v$ can be chosen linear.
\end{thm}

\vskip0.2cm

\begin{thm}\label{1:w}
Let $\BG$ be a compact graph and $\mu\in\GM(\BG)$.
\begin{itemize}
\item (i) Let $1\le p<\infty$, then 
for any $n\in\N$ there exists a linear operator 
$P_n:\tL^{1,p} (\BG)\to\Step(\BG)$ such that 
$\rank (P_n)\le n$
and 
\begin{equation}\label{1:2}
\|u-P_n u\|_{p,\mu}\le \frac{|\BG|^{1/{p'}}
\mu(\BG)^{1/p}}{n+1}
\|u'\|_p,\qquad \forall u\in\tL^{1,p} (\BG).
\end{equation}
\item(ii) Let $p=\infty$ and $d\mu=Vdx$ where $V\in\tL^\infty(\BG)$. Then
for any $n\in\N$ there exists a linear operator 
$P_n:\tL^{1,\infty} (\BG)\to\Step(\BG)$ such that 
$\rank (P_n)\le n$
and 
\begin{equation}\label{1:2a}
\|u-P_n u\|_{\infty,V}\le \frac{|\BG|\|V\|_\infty}{n+1}
\|u'\|_\infty,\qquad \forall u\in\tL^{1,\infty} (\BG).
\end{equation}
\end{itemize}
\end{thm}
In Subsection 6.1 we show in particular that the factor 
$(n+1)^{-1}$ in \eqref{1:2}
and \eqref{1:2a}
is the best possible for each $n$. 
\vskip0.2cm

The simplest example of a metric graph is the single segment $[0,L]\subset\R$.
For this case, above theorems basically turn into the results of Theorems
3.1 and 3.3 of the paper 
\cite{BS1} by M. Sh. Birman and the author (more exactly, 
into the one-dimensional particular case of these results). The most 
important feature of the estimates \eqref{1:1}, \eqref{1:2} and \eqref{1:2a} 
is their uniformity with respect to all graphs of a given length.

Our proofs are based upon \thmref{2:par} on partitioning
of a graph. This
theorem can be considered as a far going generalization of Theorem 4.1
from \cite{BS1}.  For trees and absolutely continuous
measures $d\mu=Vdx$ \thmref{2:par} 
was established in \cite{S}.

\vskip0.2cm
Let us describe the structure of the paper. The auxiliary result about 
partitioning of graphs is stated in Section 2, its proof is postponed until
Section 5. In Section 3 we prove Theorems \ref{1:uni} and \ref{1:w}, more 
exactly we are dealing with their generalizations to the Sobolev spaces
with weights. In Section 4 we consider Sobolev spaces of fractional order
and prove the corresponding analogs of Theorems \ref{1:uni} and \ref{1:w}.

The final Section 6 is devoted to discussion of the results obtained. 
In particular, we interpret 
our results in terms of approximation numbers of the appropriate embedding 
operators. We also show that in the case when $\BG$ is a tree
\thmref{1:w} and its generalization, \thmref{3:ww}, can be translated into 
the language of Hardy-type integral operators. The behavior of 
approximation numbers of such operators was studied in detail in \cite{EHL},
and there are some important intersections between our corresponding
results. We discuss them in Subsection 6.5.

For $p=2$, the results about approximation can be reformulated in terms of
the eigenvalue estimates for certain compact operators in a Hilbert space. 
In the present paper we do not
touch upon this problem. For the most important case
of \thmref{1:w} and absolutely continuous measures $\mu$ this was done 
in \cite{S}, and similar applications of our other results 
can be obtained in the same way.

\vskip0.2cm

The author  expresses his thanks
to Professor W. D. Evans for the fruitful discussions.

\section{ The key auxiliary result }
 
Let $\BG$ be a compact graph. We always consider connected graphs, 
including the ones with loops and multiple joins. For two
vertices $v,w$ the notation $v\sim w$ means that there exists an edge
$e\in\CE$ whose ends are $v$ and $w$. 
Connectedness of the graph means that for any two vertices 
$v,w\in\CV$, $v\neq w$ there exists a sequence 
$\{v_k\}_{0\le k\le m}$ of 
vertices, such that $v_0=v,\ v_m=w$ and $v_{k-1}\sim v_k$ for each 
$k=1,\ldots,m$. The {\sl combinatorial distance} $\rho_{\comb}(v,w)$ is
defined as the minimal possible $m$ in this construction. We let
$\rho_{\comb}(v,v)=0$ for any $v\in\CV$. 

The {\sl degree} $d(v)$ of a vertex $v$ is the total number
of edges incident to $v$. 
The graphs $\BG$ consisting of a single vertex 
(i.e. $\#\CV(\BG)=1,\
\CE(\BG)=\emptyset$) are called degenerate. If the (connected) graph 
$\BG$ is non-degenerate,
then its vertices $v$ with $d(v)=1$ form its
boundary $\p\BG$.

We say that a graph $G$ is a {\sl subgraph of $\BG$} 
if $G$ is a closed and connected subset of $\BG$.
According to this definition, the vertices of a subgraph
not necessarily are vertices of the original graph. For this reason,
it is often convenient to treat an arbitrary point $x\in\BG$ as a vertex.
We set $d(x)=2$ for any $x\notin\CV(\BG)$ and write $v\sim x$ if $v\in\CV(\BG)$
is one of the endpoints of the edge containing $x$. 
Given a subgraph $G$, we denote by $d_G(x)$ the degree of a
point $x\in G$ with respect to $G$. Clearly, always $d_G(x)\le d(x)$.
Note also that $\rho_G(x,y)\ge\rho_\BG(x,y)$ for any $x,y\in G$.

Along with subgraphs, our constructions involve
arbitrary connected, not necessarily closed subsets $E\subset\BG$.
Below $\CC(\BG)$ stands for the set of all such subsets. If
$E\in\CC(\BG)$, then the closure $\overline{E}$ is a subgraph, and the
complement $\overline{E}\setminus E$ is a finite set. The distinction
between $E$ and $\overline{E}$ is important only when dealing with
measures $\mu\in \GM(\BG)$ having non-zero point charges.
\vskip0.2cm

We denote by $\sqcup$ the union of subsets which are mutually disjoint, and
say that the subsets $E_1,\ldots,E_k\in \CC(\BG)$ form {\sl a partition}, or 
 {\sl a splitting} of a set $E\in \CC(\BG)$, if 
$E=E_1\sqcup\ldots\sqcup E_k$. If $E,E_1\in \CC(\BG)$ and $E_1\subset E$,
then sets $E_2,\ldots,E_k\in \CC(\BG)$ can be always found which together with 
$E_1$ form
a partition of $E$.
\vskip0.2cm

Let $\Phi$ be a non-negative function defined on the set $\CC(\BG)$
and taking values in $[0,\infty)$.
We call the function $\Phi$ {\sl super-additive} if
\begin{equation}\label{2:super}
E=\bigsqcup\limits_{j=1}^k E_j\Longrightarrow \sum\limits_{j=1}^k\Phi(E_j)
\le\Phi(E).
\end{equation} 
It is clear that any super-additive function
is monotone:  
\begin{equation}\label{2:monot}
E_1\subset E\Longrightarrow
\Phi(E_1)\le\Phi(E). 
\end{equation}

We are interested in the class $\tS(\BG)$ consisting of all
super-additive functions satisfying some additional properties which are
listed below.

 {\bf 1)} Let $\{E^r\},\ r\in\N$ be a family of sets from $\CC(\BG)$.
Then
\begin{gather} 
\Phi(E^r)\to\Phi\bigl(\cap_n E^n\bigr)\ {\text{as}}\ r\to\infty
\qquad{\text{if}}\ E^1\supset E^2\supset\ldots;\label{2:sup}\\
\Phi(E^r)\to\Phi\bigl(\cup_n E^n\bigr)\ {\text{as}}\ r\to\infty
\qquad{\text{if}}\ E^1\subset E^2
\subset\ldots.\label{2:sub}
\end{gather}

{\bf 2)} $\Phi(\{x\})=0$ for any $x\in\BG$.
\vskip0.2cm
Let $\GM_0(\BG)$ stand for the set of all measures $\mu\in\GM(\BG)$,
such that $\mu$ has no points of positive measure. It is clear that
$\GM_0(\BG)\subset\tS(\BG)$.
A more general example is given by the implication
\begin{equation}\label{2:gen}
\Phi(E)=\mu_1(E)^\a\mu_2(E)^{1-\a},\qquad 
\mu_1\in\GM_0(\BG),\ \mu_2\in\GM(\BG),\ 0<\a<1\ \Longrightarrow\ 
\Phi\in\tS(\BG).
\end{equation}
Indeed, the super-additivity of $\Phi$ is implied by H\"older's
inequality, {\bf 1)} follows from the standard properties of measures,
and {\bf 2)} follows from the condition $\mu_1\in\GM_0(\BG)$.

It is important for the applications 
that only one of two measures $\mu_1,\ \mu_2$ has to
belong to the set $\GM_0(\BG)$.
\vskip0.2cm
\vskip0.2cm
Along with partitions, we shall use {\sl pseudo-partitions}. Let $E,\G_1,
\ldots,\G_r\in\CC(\BG)$ and $E=\cup_{j=1}^r \G_j$. We say that this
is a pseudo-partition of $E$, if $\#(\G_i\cap\G_j)<\infty$ for any 
$i,j=1,\ldots,r,\ i\neq j$. We call a pseudo-partition {\sl nice} 
if the intersection $\cap_{j=1}^r \G_j$ is not empty. This intersection is
necessarily finite. 

With each function $\Phi\in\tS(\BG)$ we associate another function
$\wt\Phi$ which is defined as follows: 
\begin{equation}\label{2:til}
\wt\Phi(E)=\inf\max\limits_{j=1,\ldots,r}\Phi(\G_j)
\end{equation}
where the infimum is taken over the set of all nice pseudo-partitions 
of the set
$E$.

All our results on approximation will be derived
from the following Theorem 2.1 on super-additive functions on $\CC(\BG)$. 
\vskip0.2cm
\begin{thm}\label{2:par} 
Let $\BG$ be a compact metric graph and
$\Phi\in\tS(\BG)$. Then for any $n\in\N$ there exists a partition
$\BG=E_1\sqcup\ldots\sqcup E_k$ of $\BG$ into a family of subsets 
from $\CC(\BG)$
such that $k\le n$ and 
\begin{equation}\label{2:x}
\wt\Phi(E_j)\le (n+1)^{-1}\Phi(\BG),\qquad\forall j=1,\ldots,k.
\end{equation}
\end{thm}

The proof is rather complicated and we postpone it until Section 5.
For super-additive functions $\Phi$ such that 
\begin{equation*}
\bigl\{E,E_0\in\CC(\BG),\ |E\setminus E_0|+|E_0\setminus E|\to 0\bigr\}\
\Longrightarrow\ \bigl\{\Phi(E)\to\Phi(E_0)\bigr\}
\end{equation*}
both the formulation and the proof become much more transparent. This 
happens due to the fact that then $\Phi(E)=\Phi(\overline E)$ for any
$E\in\CC(\BG)$, and
the difference between partitions and pseudo-partitions becomes unimportant.
This simplified version of \thmref{2:par} was
obtained in \cite{S}. The general result we give here, is necessary only
for handling measures $\mu\notin\GM_0(\BG)$ in \thmref{1:w} and its
generalizations, Theorems \ref{3:ww} and \ref{5:wfr}.
\vskip0.2cm
 Now we turn to applications of \thmref{2:par}.

\section{Approximation of weighted Sobolev spaces}
\subsection{Weighted Sobolev spaces.}
Theorems \ref{1:uni} and \ref{1:w} are particular cases of similar results
for the weighted Sobolev spaces. For this reason we do not present separate 
proofs
of the original theorems but do this for the corresponding general results. We
start with the necessary definitions.

Let $\BG$ be a compact metric graph, 
$1\le p\le\infty$, and $p'=p(p-1)^{-1}$. Let $a(x)$
be a measurable function on $\BG$ such that $a(x)>0$ a.e. It is convenient to 
associate with $a(x)$ another function,
\begin{equation}\label{3:00}
w_a(x)=a(x)^{-1/p},\ p<\infty;\qquad w_a(x)=a(x)^{-1},\ p=\infty.
\end{equation}
Our basic assumption is 
$w_a\in \tL^{p'}(\BG)$.
For $p<\infty$ this is equivalent to $1/a\in \tL^{p'-1}(\BG)$.
A function $u$ on $\BG$ belongs to
the weighted Sobolev space $\tL^{1,p}(\BG,a)$ if $u$ is 
continuous on $\BG$, its restriction to each edge $e\in\CE$ 
has the distributional derivative $u'$, and 
$\|u'\|_{p,a}<\infty$. The latter functional defines on $\tL^{1,p}(\BG,a)$ a
semi-norm vanishing on the subspace $\BC$ of constant functions. It is often
convenient to factorize $\tL^{1,p}(\BG,a)$ over  $\BC$, on the
resulting quotient space $\wh{\tL}^{1,p}(\BG,a)
:=\tL^{1,p}(\BG,a)/\BC$ the 
functional $\|u'\|_{p,a}$
becomes the norm.

\subsection{Uniform approximation.}
If $a\equiv 1$, the following result  turns into \thmref{1:uni}.

\begin{thm}\label{3:uniw}
Let $\BG$ be a compact graph and let $a(x)$ be a non-negative function
on $\BG$, such that $w_a\in \tL^{p'}(\BG)$.
Then for any function
$u\in\tL^{1,p} (\BG,a)$ and any $n\in\N$ there exists a function
$v\in\Step_n(\tG)$ such that 
\begin{equation*}
\|u-v\|_\infty\le\frac{\|w_a\|_{p'}\|u'\|_{p,a}}{n+1}.
\end{equation*}
If $p=\infty$, the mapping $u\mapsto v$ can be chosen linear.

\end{thm}

\begin{proof} 1. Let first $1<p<\infty$. 
Let $L$ be a polygonal path on $\BG$ connecting 
two given points
$x_0,x$ and parametrized by the ark 
length. For any function $u\in\tL^{1,p}(\BG,a)$, 
\begin{equation*}
u(x)-u(x_0)=\int_L u'(y)dy.
\end{equation*}
Indeed, this is clearly true
if $x,x_0$ lie on the same edge, and due to the continuity of $u$ on the 
whole of $\BG$ the equality extends to any $x,x_0\in\BG$. By H\"older's 
inequality,
\begin{equation}\label{3:3}
|u(x)-u(x_0)|\le \bigl(\int_L w_a^{p'}dx\bigr)^{1/{p'}}
\bigl(\int_L a(y)|u'(y)|^pdy\bigr)^{1/p}.
\end{equation}

Given a function $u\in\tL^{1,p}(\BG,a)$,
define the function of subsets $E\in\CC(\BG)$,
\begin{equation}\label{3:4}
\Phi_u(E)=\|w_a\|_{\tL^{p'}(E)}\|u'\|_{\tL^p(E,a)}.
\end{equation}
Evidently  $\Phi_u\in\tS(\BG)$, and
$\Phi_u(\BG)=\|w_a\|_{p'}\|u'\|_{p,a}$.
It follows from \eqref{3:3} that 
\begin{equation}\label{3:5}
\sup\limits_{x\in E}|u(x)-u(x_0)|\le \Phi_u(E),\qquad\forall 
x_0\in \overline E,
\end{equation}
for any set $E\in\CC(\BG)$.

Let now $E=\G_1\cup\ldots\cup \G_r$ be
a nice pseudo-partition of $E$. According to the definition, there is a
point $x_0\in\cap_{j=1}^r \G_j$.
Applying the inequality \eqref{3:5} to each 
$\G_j$, we come to the inequality
\begin{equation*}
\sup\limits_{x\in E}|u(x)-u(x_0)|=\max\limits_{j=1,\ldots,r}
\sup\limits_{x\in \G_j}|u(x)-u(x_0)|\le \max\limits_{j=1,\ldots,r}\Phi_u(\G_j).
\end{equation*}
Minimizing the right-hand side over the set of all points $x_0\in\cap_j\G_j$
and then over the set of all nice 
pseudo-partitions of $E$
and taking into account the definition \eqref{2:til}, we find a point 
$x_E\in \overline E$
such that
\begin{equation}\label{3:61}
\sup\limits_{x\in E}|u(x)-u(x_E)|\le \wt\Phi_u(E).
\end{equation}

Suppose that the graph $\BG$ is split into the union of subsets
$E_1,\ldots,E_k\in\CC(\BG)$.
Consider the step function $v=\sum\limits_{1\le j\le k} u(x_{E_j})\chi_j$ 
where $\chi_j$
stands for the characteristic function of the set $E_j$. Then 
$v\in\Step_n(\BG)$ and by \eqref{3:61}
\begin{equation*}
\|u-v\|_\infty\le\max\limits_{j=1,\ldots,k}\wt\Phi_u(E_j).
\end{equation*}
Using \thmref{2:par}, we find a partition with $k\le n$ such that
$\wt\Phi_u(E_j)\le(n+1)^{-1}\Phi_u(\BG)$ for each $j=1,\ldots,k$. This gives 
the desired result for $1<p<\infty$.

The same argument, with minor changes, goes through for $p=1$; we skip it.
\vskip0.2cm 
2. Let now $p=\infty$, then we have instead of \eqref{3:3}:

\begin{equation*}
|u(x)-u(x_0)|\le \|au'\|_{\tL^\infty(L)}\;\int_L w_a dx \le
\|au'\|_{\tL^\infty(\BG)}\;\int_L w_a dx .
\end{equation*}
The above argument works if instead of \eqref{3:4} we take 
\begin{equation*}
\Phi(E)=\|u'\|_{\tL^\infty(\BG,a)}\int_E w_a dx.
\end{equation*}
This function of subgraphs depends on $a(x)$ but does not depend on the choice
of the function $u$. Therefore, also the partition 
$\BG=E_1\sqcup\ldots\sqcup E_k$
constructed according to \thmref{2:par} does not depend on $u$, and hence
the mapping $u\mapsto v$ is linear.
\end{proof}

\bigskip

\subsection{Weighted $\tL^p$-approximation.} Now we turn to a generalization of
\thmref{1:w}.

\begin{thm}\label{3:ww}
Let $\BG$ be a compact graph and let $a(x)$ be a non-negative function on $\BG$
such that $w_a\in\tL^{p'}(\BG)$.
\begin{itemize}
\item (i)
Let $1\le p<\infty$ and $\mu\in \GM(\BG)$. Then
for any $n\in\N$ there exists a linear operator 
$P_n=P_{n,\mu}:\tL^{1,p} (\BG,a)
\to \Step(\BG)$ such that $\rank(P_n)\le n$ and 
\begin{equation}\label{3:8}
\|u-P_n u\|_{p,\mu}\le \frac{\|w_a\|_{p'}
\mu(\BG)^{1/p}}{n+1}\|u'\|_{p,a},\qquad \forall u\in\tL^{1,p} (\BG,a).
\end{equation}

\item (ii) Let $p=\infty$ and $d\mu=Vdx$ where $V\in\tL^\infty(\BG)$. 
Then for any $n\in\N$ 
there exists a linear operator 
$P_n=P_{n,a}:\tL^{1,\infty} (\BG,a)\to \Step(\BG)$ such that $\rank(P_n)\le n$ 
and
\begin{equation*}
\|u-P_n u\|_{\infty,V}\le \frac{\|w_a\|_1
\|V\|_\infty}{n+1}\|u'\|_{\infty,a},\qquad \forall u\in\tL^{1,\infty} (\BG,a).
\end{equation*}

\end{itemize}
\end{thm} 
\begin{proof} (i) Let $1<p<\infty$; we do not discuss minor changes
needed in the case $p=1$.
The proof is quite similar to the previous one. This time we use the function
\begin{equation*}
\Phi_\mu(E)=\|w_a\|_{\tL^{p'}(E)}\mu(E)^{1/p},
\qquad E\in\CC(\BG),
\end{equation*}
cf. \eqref{3:4}. By \eqref{2:gen}, this function also lies in $\tS(\BG)$, and
$\Phi_\mu(\BG)=\|w_a\|_{p'}\mu(\BG)^{1/p}$.

Let $E=\G_1\cup\ldots\cup \G_r$ be
a nice pseudo-partition of a given subset $E\in\CC(\BG)$ and let
$x_0\in\cap_{j=1}^r \G_j$. Then we find, using \eqref{3:3}:
\begin{gather*}
\int_E|u(x)-u(x_0)|^p d\mu(x)\le\sum\limits_{j=1}^r\sup\limits_{x\in\G_j}
|u(x)-u(x_0)|^p\mu(\G_j)\\
\le \sum\limits_{j=1}^r
\bigl(\int_{\G_j}w_a^{p'}dx\bigr)^{p-1}\mu(\G_j)
\int_{\G_j} a(y)|u'(y)|^pdy \le
\biggl(\max\limits_{j=1,\ldots,r}\Phi_\mu(\G_j)\biggr)^p
\int_E a(y)|u'(y)|^pdy. 
\end{gather*}
Minimizing over the set of all points $x_0\in\cap_j\G_j$ and then 
over the set of all
nice pseudo-partitions of $E$, we find a point 
$x_{E,\mu}\in \overline E$ such that
\begin{equation}\label{3:10}
\int_E|u(x)-u(x_{E,\mu})|^p d\mu(x)\le\bigl(\wt\Phi_\mu(E)\bigr)^p
\int_E a(y)|u'(y)|^pdy.
\end{equation}
Suppose now that the graph $\BG$ is split into the union of subsets
$E_1,\ldots,E_k\in\CC(\BG)$ and let $v$ be the step function 
$v=\sum\limits_{1\le j\le k}u(x_{E_j,\mu}) \chi_j$. Then we derive from 
\eqref{3:10} that
\begin{gather*}
\int_{\BG}|u(x)-v(x)|^p d\mu(x)\le\sum_{j=1}^k \bigl(\wt\Phi_\mu(E_j)\bigr)^p
\int_{E_j} a(y)|u'(y)|^pdy \\
\le \bigl(\max\limits_{j=1,\ldots,k}
(\wt\Phi_\mu(E_j)\bigr)^p\int_{\BG} a(y)|u'(y)|^pdy.
\end{gather*}
Applying \thmref{2:par} to the function $\Phi_\mu$, we find a partition
with $k\le n$, for which 
\begin{equation*}
\max\limits_{j=1,\ldots,k}
\wt\Phi_\mu(E_j)\le (n+1)^{-1}\Phi_\mu(\BG). 
\end{equation*}
This partition depends on $\mu$ but
does not depend
on the choice of the function $u\in\tL^p(\BG,a)$. This implies that the
operator $P_n:u\mapsto v$ is linear, and we arrive at \eqref{3:8}.
\vskip0.2cm
(ii) The result is an immediate consequence of 
\thmref{3:uniw}.
\end{proof}

\section{Approximation of Sobolev spaces of fractional order} 
\subsection{Spaces $\tL^{\t,p}(\BG)$.} As in the previous Section,
it is convenient for us to consider the spaces factorized over the
subspace $\BC$ of constant functions. However, in our notations
we do not distinguish between a function $u$ and the corresponding
factor-element. In order to simplify our reasonings, we consider only
$1<p<\infty$ and the spaces without weights.

The most natural approach to the spaces $\tL^{\t,p}(\BG)$ uses interpolation
between the space $\wh\tL^{1,p}(\BG)$, see Subsection 3.1, and
the quotient space $\wh\tL^p(\BG)=\tL^p(\BG)/\BC$. 
As usual, the norm in $\wh\tL^p(\BG)$ is defined by 
\begin{equation*}
\|u\|_{\wh\tL^p(\BG)}=\min_{c\in\C}\|u-c\|_p.
\end{equation*}
The spaces $\wh\tL^p(\BG)$ and $\wh\tL^{1,p}(\BG)$ form a Banach couple,
see e.g. \cite{T}, and we define  the interpolation space
\begin{equation}\label{4:interp}
\wh\tL^{\t,p}(\BG)=\bigl(\wh\tL^p(\BG),\wh\tL^{1,p}(\BG)\bigr)_{\t,p},
\qquad 0<\t<1.
\end{equation}
We write $u\in \tL^{\t,p}(\BG)$, when it is convenient to view $u$
as an individual function rather than the equivalence class $\{u+\BC\}$.
We do not discuss here interpolation with the second parameter
$q\neq p$ which would lead to the general Besov spaces.

There are many ways to define an interpolation norm in $\wh\tL^{\t,p}$. 
For our purposes
it is convenient to
use the $L$-method with the parameters $p_0=p_1=p$, see e.g.
\cite{T}, Section 1.4. So, we define for $0<t<\infty$:
\begin{equation}\label{4:l}
L(t,u;\BG)=\inf\bigl\{\|u_0\|_{\wh\tL^p(\BG)}^p+t\|u'_1\|_{\tL^p(\BG)}^p:
\ u=u_0+u_1;\
u_0\in\tL^p(\BG), u_1\in\tL^{1,p}(\BG)\bigr\}.
\end{equation}
A function $u\in\tL^p(\BG)+\tL^{1,p}(\BG)$ belongs to the space 
$\tL^{\t,p}(\BG)$ if and only if
\begin{equation}\label{4:sobfr}
\bigl(\|u\|_{\tL^{\t,p}(\BG)}\bigr)^p:=\int_0^\infty t^{-1-\t}L(t,u;\BG)dt
<\infty.
\end{equation}
\vskip0.2cm

Replacing in \eqref{4:sobfr} the graph $\BG$ by its arbitrary subset
$E\in\CC(\BG)$ and fixing an element $u\in\tL^{\t,p}(\BG)$, we obtain the 
function 
\begin{equation}\label{4:10}
J_{\t,u}(E)=\bigl(\|u\|_{\wh\tL^{\t,p}(E)}\bigr)^p,\qquad E\in\CC(\BG). 
\end{equation}
Let us show that $J_{\t,u}\in\tS(\BG)$. 
First of all, we note that
the function $
J_{0,u}(E)=\|u\|_{\wh\tL^p(E)}^p $
lies in $\tS(\BG)$.  
Indeed, the properties {\bf 1)} and {\bf 2)}
of functions $\Phi\in\tS(\BG)$, cf. Section 2, are evidently
fulfilled, and
for any constant $c$ and any partition $E=E_1\sqcup\ldots\sqcup E_k$ we have
\begin{equation*}
\int_E|u-c|^pdx=\sum_{j=1}^k\int_{E_j}|u-c|^pdx\ge
\sum_{j=1}^k\inf_{c_j\in\C}\int_{E_j}|u-c_j|^pdx
\end{equation*}
which yields super-additivity. The function $J_{1,u}(E)=\|u'\|_{\tL^p(E)}^p$
also lies in $\tS(\BG)$, therefore the same is true for the function
$L(t,u;G)$ defined by the equality \eqref{4:l} for the set $E\in\CC(\BG)$
substituted for $\BG$. 
Integration in \eqref{4:sobfr} does not violate the
property of a function to lie in $\tS(\BG)$. Hence, it is proved
that 
$J_{\t,u}\in\tS(\BG)$.
It follows from here and \eqref{2:monot} that 
\begin{equation}\label{4:b}
\|u\|^p_{\wh\tL^{\t,p}(E_1)}\le\|u\|^p_{\wh\tL^{\t,p}(E)},\qquad
\forall E,E_1\in\CC(\BG),\ E_1\subset E.
\end{equation}
\vskip0.2cm
If $\t p>1$, any function 
$u\in \tL^{\t,p}(\BG)$ is continuous. This is well known when 
$\BG$ is
a single segment. Hence, $u$ is continuous
on any polygonal path in $\BG$
and thus, on the whole of $\BG$. 

Denote by $C(\t,p)$ the sharp constant in the inequality
\begin{equation}\label{4:2}
\max\limits_{x,x_0\in[0,l]}|u(x)-u(x_0)|\le C(\t,p)l^{\t-1/p}
\|u\|_{\tL^{\t,p}[0,l]}.
\end{equation}
The value of $C(\t,p)$ does not depend on $l$, which follows from the
homogeneity arguments. The inequality \eqref{4:2} automatically extends to
the graphs: due to \eqref{4:b}, 
\begin{equation}\label{4:3}
\sup\limits_{x,x_0\in E}|u(x)-u(x_0)|\le C(\t,p)|E|^{\t-1/p}J_{\t,u}(E)^{1/p},
\qquad\forall\ E\in\CC( \BG)
\end{equation}
where the function $J_{\t,u}(E)$ is defined by \eqref{4:10}.
\subsection{Approximation of $\tL^{\t,p}$.}
Below are analogs of Theorems \ref{1:uni} and \ref{1:w} for 
the spaces $\tL^{\t,p}(\BG)$.

\begin{thm}\label{5:unifr}
Let $\BG$ be a compact graph, $0<\t<1$, and $1/\t<p<\infty$.  
Then for any function
$u\in\tL^{\t,p} (\BG)$ and any $n\in\N$ there exists a function
$v\in\Step_n(\BG)$ such that
\begin{equation}\label{5:11}
\|u-v\|_\infty\le C(\t,p)\frac{|\BG|^{\t-1/p}\|u\|_{\tL^{\t,p}(\BG)}}
{(n+1)^\t}.
\end{equation}
\end{thm}

We only outline the proof; details can be easily reconstructed
by analogy with \thmref{3:uniw}.

Together with $J_{\t,u}(E)$, the function 
\begin{equation}\label{5:110}
\Phi_u(E)=|E|^{1-1/{(p\t)}}J_{\t,u}(E)^{1/{(p\t)}}
\end{equation}
also belongs to $\tS(\BG)$,
cf. \eqref{2:gen}. Let subsets $\G_1,\ldots,\G_r\in\CC(\BG)$ form a nice
pseudo-partition of a set $E\in\CC(\BG)$ and let $x_0$ be a point from 
their intersection.
The inequality 
\begin{equation*}
\sup\limits_{x\in E}|u(x)-u(x_0)|\le C(\t,p)\bigl(\max\limits_{j=1,\ldots,r}
\Phi_u (\G_j)\bigr)^\t
\end{equation*}
is easily derived from \eqref{4:3}.
Minimizing over the set of all points from $\cap_j\G_j$ and
then, over the set of all nice pseudo-partitions of $E$,
we find a point $x_E\in \overline E$, such that
\begin{equation*}
\sup_{x\in E}|u(x)-u(x_E)|\le C(\t,p)\bigl(\wt\Phi_u (E)\bigr)^\t.
\end{equation*}
The proof is concluded by 
applying \thmref{2:par} to the function \eqref{5:110}.
\qed
\vskip0.2cm

\begin{thm}\label{5:wfr}
Let $\BG$ be a compact graph, $0<\t<1$, and $1/\t<p<\infty$. Let 
$\mu\in\GM(\BG)$.
Then for any $n\in\N$ there exists a linear operator 
$P_n:\tL^{\t,p} (\BG)\to\Step(\BG)$ such that 
$\rank (P_n)\le n$
and 
\begin{equation}\label{5:21}
\|u-P_n u\|_{p,\mu}\le C(\t,p)\frac{|\BG|^{\t-1/p}
\mu(\BG)^{1/p}}{(n+1)^\t}
\|u\|_{\tL^{\t,p}},\qquad \forall u\in\tL^{\t,p} (\BG).
\end{equation}
\end{thm}

Again, we only sketch the proof.
We make use of the function
\begin{equation*}
\Phi_\mu(E)=|E|^{1-1/(\t p)}\mu(E)^{1/(\t p)}
\end{equation*}
which by \eqref{2:gen} belongs to $\tS(\BG)$. For any set $E\in\CC(\BG)$
we find a point $x_{E,\mu}\in \overline E$ such that 
\begin{equation}\label{5:48}
\int_E|u(x)-u(x_{E,\mu})|^p d\mu(x)\le C(\t,p)^p
\bigl(\wt\Phi_\mu(E)\bigr)^{\t p}
J_{\t,u}(E),
\end{equation}
cf. \eqref{3:5}. 
Let $\BG=E_1\sqcup\ldots\sqcup E_k$ be an arbitrary partition of the graph
$\BG$.  Let $v=\sum_{j=1}^k u(x_{E_j,\mu})\chi_j$, then we derive from 
\eqref{5:48}  
using the super-additivity of $J_{\t,u}$: 
\begin{gather*}
\int_\BG|u-v|^p d\mu(x)=\sum_{j=1}^k\int_{E_j}|u-u(x_{E_j,\mu})|^p d\mu(x)
\le C(\t,p)^p \bigl(\wt\Phi_\mu(E_j)\bigr)^{\t p}J_{\t,u}(E_j)\\ 
\le C(\t,p)^p 
\biggl(\max\limits_{j=1,\ldots,k}\wt\Phi_\mu(E_j)\biggr)^{\t p}J_{\t,u}(\BG).
\end{gather*}
We come to the desired result applying \thmref{2:par} to the 
function $\Phi_\mu$
and taking into account that the mapping $P:u\mapsto v$ is linear.
\qed

\section{Proof of \thmref{2:par}}

\subsection{The case of trees.} Let $\BG=\BT$ be a tree, that is  
connected graph without cycles, loops and multiple joins. 
For any two points $x,y\in\BT$ there exists a unique simple
polygonal path  in $\BT$ connecting $x$ with $y$, we denote it by 
$\lu x, y\ru$. It is clear
that $|\lu x,y\ru|=\rho(x,y)$.

For trees the notion of nice pseudo-partition simplifies. 
Indeed, if $T=\T_1\cup\ldots\cup\T_r$ is a nice pseudo-partition
of a (closed) subtree $T\subset\BT$, then the intersection 
$\Xi=\cap_{j=1}^r\T_j$
consists of exactly one point. For if $x_1\neq x_2$ and
$x_1,x_2\in \Xi$, then
also $\lu x_1, x_2\ru\subset\Xi$ which contradicts the definition of 
pseudo-partition. So, the point $x\in\cap_j\T_j$ 
is uniquely defined by a nice pseudo-partition. Besides, all the subsets
$\T_j$ are necessarily closed, i.e. each of them is a subtree of $T$.

Conversely, each nice pseudo-partition of $T$ is uniquely determined by 
the choice of
the point $x$. Indeed, the tree $T$ splits in a
unique way into the union of subtrees $\T_j\subset T$, $j=1,\ldots,
d_T(x)$, rooted at $x$ and such that $d_{\T_j}(x)=1$ for each $j$. Evidently
this pseudo-partition is nice. We call
the pair $\{T,x\}$ a {\sl punctured subtree} and the above constructed
partition -- its {\sl canonical pseudo-partition}. 

Let $\Phi\in\tS(\BT)$. Defining
\begin{equation}\label{8:can}
\Phi^{\bpr}(T,x)=\max\limits_{j=1,\ldots,d_T(x)}\Phi(\T_j),
\end{equation}
we evidently have

\begin{equation}\label{8:1}
\wt\Phi(T)=\min\limits_{x\in T} \Phi^{\bpr}(T,x).
\end{equation}

Let in particular $T=\BT$. Each subtree $\T_j$ appearing in the canonical
pseudo-partition of $\{\BT,x\}$ is determined by indication of its initial
edge $\lu x, v\ru$, $v\sim x$ and we denote this subtree by $\T_{\lu
x,v\ru}$. For $T=\BT$ the definition \eqref{8:can} takes the form
\begin{equation*}
\Phi^{\bpr}(\BT,x)=\max_{v\sim x}\Phi(\T_{\lu x,v\ru}).  
\end{equation*}
 
The following lemma is the heart of our proof of \thmref{2:par}.
\begin{lem}\label{8:lem} 
Let $\BT$ be a compact
metric tree and $\Phi\in\tS(\BT)$. Then for any $\vare\in(0,\Phi(\BT))$
there exists a pseudo-partition $\BT=T\cup T'$, such that the set
$T'\setminus T$ is connected (that is, belongs to $\CC(\BT)$) and
for the single point $x^*\in T\cap T'$ the inequalities hold:
\begin{gather}
\Phi^{\bpr}(T,x^*)\le\vare;\label{8:04} \\
\Phi(T'\setminus\{x^*\})\le \Phi(\BT)-\vare. \label{8:03}
\end{gather}

\end{lem} 
\begin{proof} 
Without loss of generality, we can assume
$\Phi(\BT)=1$. Take any vertex $v_0\in\p\BT$, then
$\Phi^{\bpr}(\BT,v_0)=\Phi(\BT)=1$. There is a unique vertex $v_1\sim v_0$.
Now we choose the vertices $v_2\sim v_1,\ldots, v_{k+1}\sim v_k,\ldots$ as
follows. If $v_k$ is already chosen, we define $v_{k+1}$ as the vertex
different from $v_{k-1}$ and such that
 \begin{equation}\label{8:path} 
\Phi(\T_{\lu v_k,v_{k+1}\ru})= \max_{w\sim
v_k,w\neq v_{k-1}}\Phi(\T_{\lu v_k,w\ru}) =\Phi^{\bpr}(\T_{\lu
v_k,v_{k+1}\ru},v_k). 
\end{equation} 
If there are several vertices $w\sim v_k$ at which the maximum in 
the middle term of \eqref{8:path} is
attained, then any of them can be chosen as $v_{k+1}$. The described
procedure is always finite, it terminates when we arrive at a vertex
$v_m\in\p\BT$. On the path $\CP= \lu v_0,v_m\ru$ we
introduce the natural ordering, i.e. $y\succeq x$ means that $x\in\langle
v_0,y\rangle$. We write $y\succ x$ if $y\succeq x$ and $y\neq x$. 
 
Let $x\in \CP$ be not a vertex of $\BT$, then $v_{k-1}\prec x\prec
v_k$ for some $k=1,\ldots,m$. Denote 
\begin{equation*} 
T_x^+=\T_{\lu x,v_k\ru},\qquad T_x^-=\T_{\lu x,v_{k-1}\ru},\qquad 
x\neq v_0,\ldots,v_m.
\end{equation*}
We also define the subtrees $T_x^\pm$ for $x=v_0,\ldots, v_m$. Namely,
\begin{gather*} 
T_{v_k}^-=T_{\langle v_k,v_{k-1}\rangle},\qquad k=1,\ldots,
m;\\ T_{v_0}^+=\BT,\qquad
T_{v_k}^+=\bigcap_{v_{k-1}\prec x\prec v_k}T_x^+= \bigcup_{v\sim v_k,
v\neq v_{k-1}}T_{\langle v_k,v\rangle}, \qquad k=1,\ldots, m-1.
\end{gather*} 
Finally, $T_{v_0}^-=\{v_0\},\ T_{v_m}^+=\{v_m\}$ are
degenerate subtrees. For any $x\in\CP$ we have
$\BT=T_x^+\cup T_x^-$. Clearly, this 
is a pseudo-partition of the tree $\BT$, and $T_x^+\cap
T_x^-=\{x\}$. Besides, for any $x\in\CP$ we have $x\in\p T_x^-$, and the set  
$T_x^-\setminus T_x^+=T_x^-\setminus\{x\}$ is connected, i.e.
belongs to $\CC(\BT)$.
\vskip0.2cm
 
The function $F(x)=\Phi(T_x^+)$ is well defined on $\CP$ and non-increasing.
By \eqref{2:sup}, $F$ is left-continuous with 
respect to the ordering adopted. By the construction, 
\begin{equation*} 
\Phi^{\bpr}(T^+_{x_0},x_0)=F(x_0),\qquad \forall x_0\in\CP.
\end{equation*} 
 Further, \eqref{2:sub}
implies that 
\begin{equation*}
F(x_0+):=\lim_{x\succ x_0,x\to x_0}F(x)=\Phi(T_{x_0}^-\setminus\{x_0\}),
\qquad \forall x_0\in\CP.
\end{equation*} 
We also have 
\begin{equation*}
0=F(v_m)<\vare<F(v_0)=1. 
\end{equation*} 
Therefore, there exists a point $x^*\in \CP$ such that
\begin{equation*} 
\Phi^{\bpr}(T_{x^*},x^*)=F(x^*)\ge\vare \ge F(x^*+).
\end{equation*}
We take $T=T_{x^*}^+$ and $T'=T_{x^*}^-$. Then the inequality \eqref{8:04} 
is satisfied
and \eqref{8:03} is implied by super-additivity: 
\begin{equation*}
\Phi(T'\setminus\{x^*\})\le 1-\Phi(T)=1-F(x^*)\le 1-\vare. 
  \end{equation*} 
\end{proof}
 
\subsection { Proof of \thmref{2:par} for the case of trees.} 
Let $\BG=\BT$ be a tree.

1. Let $n=1$. Apply the
result of Lemma \ref{8:lem} with $\vare=\Phi(\BT)/2$. Let $\BT=T\cup T'$
be the corresponding pseudo-partition, then
$\Phi^{\bpr}(T',x^*)\le\Phi(T')\le\Phi(\BT)/2$.  Consider the canonical
pseudo-partition of the punctured tree $\{\BT,x^*\}$. Each subtree of this
pseudo-partition is contained either in $T$ or in $T'$, therefore
\begin{equation*}
\Phi^{\bpr}(\BT,x^*)\le\max\bigl(\wt\Phi(T,x^*),
\Phi^{\bpr}(T',x^*)\bigr)\le\Phi(\BT)/2. 
\end{equation*} 
Taking \eqref{8:1} into account, we see that \eqref{2:x} with $k=n=1$ is
satisfied if we take $E_1=\BT$. 
 
\vskip0.2cm 2. We proceed by induction. Suppose that the result is already
proved for $n=n_0-1$.  Let $\BT=T\cup T'$ be the pseudo-partition constructed
according to Lemma \ref{8:lem} for $\vare=(n_0+1)^{-1}\Phi(\BT)$
and let $T\cap T'=\{x^*\}$. Then
\begin{equation*} 
\Phi(T'\setminus\{x^*\})\le n_0(n_0+1)^{-1}\Phi(\BT). 
\end{equation*} 
Let us define a function $\Phi'$ of subsets $E\in\CC(T')$, taking
\begin{equation*} 
\Phi'(E)=\Phi(E\setminus\{x^*\}),\qquad \forall E\in\CC(T'),
\end{equation*}
then evidently $\Phi'\in\tS(T')$.
By the inductive hypothesis, there exists a splitting of
$T'$ into the union of subsets $E_j\in\CC(T')$,
$j=1,\ldots,k$ such that $ k\le n_0-1 $ and for each $j$ 
\begin{equation*}
\wh\Phi'(E_j)\le n_0^{-1}\Phi'(T')=n_0^{-1}\Phi(T'\setminus
\{x^*\})\le (n_0+1)^{-1}\Phi(\BT). 
\end{equation*}
The point $x^*$ lies in only one of the sets $E_j$, let it be $E_k$. 
Since $x^*\in \p T'$, we conclude that the set $E_k\setminus\{x^*\}$ is
connected and therefore belongs to $\CC(\BT)$.
 
The family $E_1,\ldots,E_{k-1},E_k\setminus\{x^*\},T$ forms 
the desired partition of $\BT$ 
for $n=n_0$. For the trees, the proof of \thmref{2:par} is complete.

\subsection{General case.} \thmref{2:par} for arbitrary graphs
can be easily reduced to the 
case of trees by means of ``cutting cycles''. Below we describe the
procedure of such reduction.

Let $\BG$ be a compact graph and $\Phi$ be a function from $\tS(\BG)$.
Let $e$ be an edge of $\BG$ which is a part of a cycle. Supposing that
$e$ is not a loop, we identify $e$ with the segment $[0,|e|]$. Take any
point $x\in\Int (e)$ and replace it by the pair $x_1,x_2$ of new 
vertices.
Respectively, the edge $e$ is replaced by the pair $e_1,e_2$ of new edges 
whose total length is equal to $|e|$. As the
result, we obtain a new graph, say $\BG_1$. Note that the edges $e_1,e_2$
are parts of no cycle in $\BG_1$.
Define the mapping
$\tau_1:\BG_1\to\BG$ which is identical on $\BG\setminus \Int (e)$
and sends isometrically $e_1$ onto $[0,x]$ and $e_2$ onto $[x,|e|]$.
The mapping $\tau_1$ is one-to-one on $\BG_1\setminus\{x_1,x_2\}$, and
$\tau_1(x_1)=\tau_1(x_2)=x$. It is clear that $\tau_1$ is
non-expanding and hence, continuous.

The changes in this construction, needed if $e$ is a loop, are evident.

Now, define a function $\Phi_1$ on the set $\CC(\BG)$, namely
\begin{equation*}
\Phi_1(E)=\Phi(\tau_1(E))\ {\text{if}}\ x_1\in E,\qquad
\Phi_1(E)=\Phi(\tau_1(E)\setminus\{x\})\ {\text{if}}\ x_1\not\in E.
\end{equation*}
The function $\Phi_1$ is super-additive. Indeed, let $E\in\CC(\BG_1)$ and
$E=\sqcup_{j=1}^k E_j$.
If $x_1\notin E$, then also $x_1\notin E_j$ for any $j$, and if $x_1\in E$,
then $x_1\in E_{j_0}$ for exactly one value of $j$. In both cases, the
inequality \eqref{2:super} for $\Phi_1$ is implied by the similar inequality 
for $\Phi$. The properties {\bf 1), 2)} for the function $\Phi_1$ also
follow from the same properties for $\Phi$. Hence,
$\Phi_1\in\tS(\BG_1)$.

Repeating this procedure, we obtain a sequence of graphs
$\BG_0:=\BG,\BG_1,\ldots,\BG_m$, a sequence of mappings $\tau_j:\BG_j
\to\BG_{j-1},\ j=1,\ldots,m$, and a family of functions $\Phi_j\in\tS(\BG_j)$.
The procedure stops as soon as we come to a
graph without cycles and loops, that is when $\BG_m=:\BT$ is a compact tree.
The mapping $\tau=\tau_m\circ\ldots\circ\tau_1:\BT\to\BG$ is continuous
and measure preserving. Due to the continuity of $\tau$, $T\in\CC(\BT)
\Longrightarrow \tau(T)\in\CC(\BG)$. Besides, $\tau$ transforms partition
into partition and preserves the property of a partition to be nice.
The function $\Phi_m$ belongs to $\tS(\BG_m)$ and $\Phi_m(\BT)=\Phi(\BG)$.
\vskip0.2cm
By the result of previous section, for a given $n\in\N$ there exists a 
partition $\BT=\sqcup_{j=1}^k E_j$ 
into the union of subsets from $\CC(\BT)$, such that
$k\le n$ and $\wt\Phi_m(E_j)\le (n+1)^{-1}\Phi_m(\BT)=(n+1)^{-1}\Phi(\BG)$ for 
each $j$. 
Taking $E'_j=\tau(E_j)$, we find a partition of $\BG$ which meets all
the requirements of \thmref{2:par}.   \qed

\section{Complements and concluding remarks}
\subsection{On the sharpness of estimates.} 

a) The factor $(n+1)^{-1}$ in the inequality \eqref{2:x} of 
\thmref{2:par}
is sharp for each $n$. To see this, consider the star graph $\BG_N$ 
consisting of $N$ edges
$e_k=\lu o,v_k\ru,\ k=1,\ldots,N$ of equal length $1$, all emanating from the 
root $o$. For any subset $E\in\CC(\BG_N)$ we define $\Phi(E)=|E|$,
then $\Phi\in\tS(\BG_N)$.
Take $n=N-1$, then at least one of the subsets $E_j$
appearing in the conclusion of \thmref{2:par} necessarily contains two 
edges of $\BG_N$. Thus, $\Phi(E_j)\ge 2$ and hence, $\wt\Phi(E_j)\ge1$ 
for any nice pseudo-partition of $E_j$. Since $|\BG_N|=N=n+1$, we see that 
the inequality \eqref{2:x} turns into equality.

b) The same factor $(n+1)^{-1}$ in the inequality \eqref{1:2} 
of \thmref{1:w} is also
sharp for each $n$. Indeed, consider the star graph $\BG_N$ and
the measure $\mu\in\GM(\BG_N)$
defined as $\mu=\d_{v_1}+\ldots+\d_{v_N}$. Consider also the subspace 
$Y\subset\tL^{1,p}(\BG_N)$ formed by 
the functions $u$ such that $u\res e_k=c_k\rho(o,x),\ k=1,\ldots,N$. Then
\begin{equation*}
\|u'\|_p=\|u\|_{\tL^p(\BG_N,\mu)}=\|c\|_{\ell^p_N},\qquad 
c=\{c_k\}_{1\le k\le N},\ \forall u\in Y.
\end{equation*}
It follows that for any linear operator $P:\tL^{1,p}(\BG)\to\tL^p(\BG_N,\mu)$
with $\rank(P)\le n$ the quantity 
\begin{equation*}
\inf\limits_{u\in\tL^{1,p}(\BG):\|u'\|_p=1}
\|u-Pu\|_{\tL^p(\BG_N,\mu)}
\end{equation*}
 is no smaller than the $n$-width in $\ell_N^p$
of the unit ball of this space. For $n<N$ this $n$-width is equal to one,
see e.g. \cite{P}, Proposition 1.3. Since $|\BG_N|=\mu(\BG_N)=N$, 
we see that for $n=N-1$
an element $u\in\tL^{1,p}(\BG):\|u'\|_p=1$ can always be found in such a way
that
\begin{equation*}
\|u-Pu\|_{\tL^p(\BG_N,\mu)}\ge1=\frac{|\BG_N|^{1/{p'}}\mu(\BG_N)^{1/p}}
{n+1}.
\end{equation*}

Replacing
the above measure $\mu$ by a sequence of 
measures $V_jdx$ which $*$-weakly approximate $\mu$, we find
that the factor $(n+1)^{-1}$ in \eqref{1:2} is the least possible
also for absolutely continuous measures. However, for each particular
absolutely continuous measure $\mu$ the inequality in \eqref{1:2} is
always strict.

c) The same factor in the inequality \eqref{1:1} is sharp for $n=1$. 
For $n>1$ it becomes sharp, provided one passes to the version of 
\thmref{1:uni} (and its generalization, \thmref{3:uniw}) dealing with 
vector-valued functions. Namely, let
$X$ be a Banach space and let $\tL^{1,p}(\BG;X)$ stand for the space of 
$X$-valued
functions on $\BG$ whose definition is clear by analogy with the case 
of scalar-valued functions, cf.
Section 1. Both mentioned theorems extend to the spaces $\tL^{1,p}(\BG;X)$, 
the proof actually remains the same. 

Now, take $X=\ell^\infty$. For $k\in\N$, let $\eta_k\in\ell^\infty$ be the 
element whose $k$-th coordinate is $1$ and all the others are equal to zero.  
On the star graph $\BG_N$ consider the function $u$ which is 
$\eta_k\rho(0,x)$ on the 
edge $e_k\in\BG_N$. Then $u\in\tL^{1,p}(\BG_N;X)$ for each $p\in[1,\infty]$ and
$\|u'\|_{\tL^p(\BG_N;X)}=1$. The same reasoning as in a) shows that for
$n=N-1$ the constant factor $(n+1)^{-1}$ in the vector-valued version of 
\eqref{1:1} is the best possible.

\subsection{Graphs and trees: comparison of the corresponding results.} 
Given a compact graph $\BG$, let $\BT$ and
$\tau:\BT\to\BG$ be the tree and the mapping constructed in Subsection 5.3. 
Let $a(x)$ be a non-negative function on $\BG$ such that 
$w_a\in\tL^{p'}(\BG)$ (cf. \eqref{3:00}). Define $b(x)=a(\tau(x))$, then
$w_b\in\tL^{p'}(\BT)$ and $\|w_b\|_{\tL^{p'}(\BT)}=\|w_a\|_{\tL^{p'}(\BG)}$.
Moreover, it is clear from the construction 
that the mapping $u(x)\mapsto v(x)=u(\tau(x))$
defines an isometry between the space $\tL^{1,p}(\BG,a)$ and an
appropriate subspace of finite codimension in $\tL^{1,p}(\BT,b)$. 
Indeed, suppose that the passage from the graph $\BG$ to the tree $\BT$
consists in replacing the points $x^{(j)}\in\BG$, $j=1,\ldots,m$  
by the pairs $\{x_1^{(j)},x_2^{(j)}\}\subset\BT$. Then the 
space $\tL^{1,p}(\BG,a)$ can be identified with the subspace
\begin{equation*}
\{u\in\tL^{1,p}(\BT,b): \ u(x_1^{(j)})=u(x_2^{(j)}),\ j=1,\ldots,m.\}
\end{equation*}
The above mapping $u\mapsto v$ defines also the
natural isometry between the spaces $\tL^p(\BG,V)$ and $\tL^p(\BT,W)$
where $W(x)=V(\tau(x))$.
It follows from these remarks
that \thmref{3:ww} for general graphs reduces to its particular case for trees.

The same is true for \thmref{5:wfr}, though for the spaces $\tL^{\t,p}$ the 
above mapping $u\mapsto v$ is not necessarily an isometry. But this is
always a contraction, so that the constant in the estimate \eqref{5:21}
for a graph $\BG$ can not exceed the one for the corresponding tree $\BT$.

\subsection{Approximation numbers of embedding operators.} Suppose that a point
$o\in\BG$ is fixed, and define the spaces
\begin{equation*}
\tW^{1,p}(\BG,a;o)=\{u\in\tL^{1,p}(\BG,a):u(o)=0\}
\end{equation*}
and, for $0<\t<1$ and $p>1/\t$,
\begin{equation*}
\tW^{\t,p}(\BG;o)=\{u\in\tL^{\t,p}(\BG):u(o)=0\}.
\end{equation*}
We take $\|u'\|_{p,a}$ as the norm in $\tW^{1,p}(\BG,a;o)$ and 
$\|u\|_{\tL^{\t,p}}$ (cf. \eqref{4:sobfr}) as the norm in $\tW^{\t,p}(\BG;o)$.
It is clear that the spaces $\tW^{1,p}(\BG,a;o)$ and $\tW^{\t,p}(\BG;o)$
are naturally isometric to the quotient spaces $\wh\tL^{1,p}(\BG,a)$ and 
$\wh\tL^{\t,p}(\BG)$ respectively. For this reason, Theorems 
of Sections 3 and 4 immediately apply to
the spaces $\tW^{1,p}(\BG,a;o)$ and $\tW^{\t,p}(\BG;o)$. 
\vskip0.2cm
Given two Banach spaces $Y$ and $X$ and an integer $n\ge0$, let $\tP_n$ stand 
for the set of all linear mappings $P:Y\to X$ whose rank
does not exceed $n$. Recall the definition of the approximation numbers 
$a_n(T)$ of a bounded 
linear operator $T:X\to Y$, 
see e.g. \cite{ET}:
\begin{equation}\label{6:3}
a_n(T)=\inf_{P\in\tP_{n-1}(Y,X)}\|T-P\|_X.
\end{equation}
In particular, this definition applies to the case when
$Y$ is embedded in $X$
algebraically and topologically, and $T=J_{Y,X}$ is the corresponding
embedding operator. \thmref{3:ww} implies that under its assumptions
we have, for any $n\in\N$:
\begin{gather}
a_n(J_{\tW^{1,p}(\BG,a;o),\tL^p(\BG,\mu)})\le \frac
{\|w_a\|_{p'}\mu(\BG)^{1/p}}{n},
\qquad p<\infty;\label{6:2}\\
a_n(J_{\tW^{1,\infty}(\BG,a;o),\tL^\infty(\BG,V)})\le\frac 
{\|w_a\|_1\|V\|_\infty}{ n}.\label{6:2a}
\end{gather}
In the same way, it follows from \thmref{5:wfr} that
\begin{equation*}
a_n(J_{\tW^{\t,p}(\BG,a;o),\tL^p(\BG,\mu)})\le C(\t,p)|\BG|^{\t-1/p}
\mu(\BG)^{1/p}n^{-\t},
\qquad \forall n\in\N,\ 1<p\t<\infty.
\end{equation*}

\subsection{Hardy-type operators on trees.} For the case of trees
there is a useful interpretation of the estimates \eqref{6:2} and \eqref{6:2a}
in terms of approximation numbers of certain integral operators.

Let $\BT$ be a compact metric tree on which a point  
$o$ (the root) is selected.  Below we use the notation $\lu x,y\ru$
introduced is Subsection 5.1. 

The Hardy-type integral operator 
with weights $v,w$ on the rooted tree $\{\BT,o\}$ is defined as
\begin{equation}\label{6:ho}
g(x)=\bigl(H_{v,w}f)(x)=
\bigl(H_{v,w}(\BT,o)f)(x)=v(x)\int_{\lu o,x\ru}f(y)w(y)dy.
\end{equation}
At first we assume that $w(x)\neq 0$ a.e.
and set $a(x)=|w(x)|^{-p}$, then $w=w_a$, cf. \eqref{3:00}. 
It is easy to see that the operator 
\begin{equation*} 
Q_w:f(x)\mapsto u(x)=\int_{\lu o,x\ru}f(y)w(y)dy
\end{equation*}
defines an isometry of the space $\tL^p(\BT)$ onto $\tL^{1,p}(\BT,a;o)$.
Besides, $\|g\|_p=\|Q_w f\|_{p,V}$ where $V=|v|^p$. This shows that
\begin{equation*}
a_n(H_{v,w})=a_n(J_{\tW^{1,p}(\BT,a;o),\tL^p(\BT,V)}),\qquad\forall n\in\N.
\end{equation*}

Now we are in a position to justify the following result.

\begin{thm}\label{6:apn}
Let $\BT$ be a compact metric tree with the root $o$
and let $w\in\tL^{p'}(\BT)$,
$v\in \tL^p(\BT)$ where $1<p<\infty$. Then the operator $H_{v,w}$ is
compact in $\tL^p(\BT)$ and its
approximation numbers satisfy the estimate

\begin{equation}\label{6:apnest}
a_n(H_{v,w})\le\frac{\|v\|_p\|w\|_{p'}}{n},\qquad \forall n\in\N.
\end{equation}
\end{thm}
\begin{proof} If $w(x)\neq 0$ a.e., then \eqref{6:apnest} immediately follows 
from \thmref{3:ww}. 
The result extends to the
general case by a standard approximation argument.
\end{proof}

\subsection{Comparison with the results of \cite{EHL}.}
The techniques of \cite{EHL}
is based upon a careful analysis of the function
$A_{v,w}(T)$ of subtrees $T\in\CC(\BT)$ which in the compact case 
can be defined
as follows:
\begin{equation*}
A_{v,w}(T)=\min_{o\in T}\|H_{v,w}(T,o):{\tL^p(T)\to\tL^p(T)}\|,
\end{equation*}
cf. Theorem 3.8 in \cite{EHL}. Evidently,
 \begin{equation*}
A_{v,w}(T)\le \|v\|_{\tL^p(T)}\|w\|_{\tL^{p'}(T)}.
\end{equation*}
Up to a change of notations, the expression in the right-hand side
is exactly the function $\Phi$ appearing in the proof of \thmref{3:ww}.
One may attempt to apply our analysis directly to the function $A_{v,w}(T)$.
However, such an attempt fails, since this function is, in general, 
not super-additive. Note also that
the converse inequality
$A_{v,w}(T)\ge c\|v\|_{\tL^p(T)}\|w\|_{\tL^{p'}(T)}$ 
with any $c>0$ is impossible.

\vskip0.2cm
In terms of the function $A_{v,w}(T)$ the authors of \cite{EHL} found 
for the approximation numbers $a_n=a_n(H_{v,w}(\BT,o))$ some two-sided
estimates, see Theorem 3.18 there. Based upon these estimates, they
justified the Weyl-type asymptotics for $a_n$.
As it was pointed out to the author 
by W.D. Evans, 
the inequality 
\begin{equation*}
a_{n+4}\le 3n^{-1}\|v\|_p\|w\|_{p'}
\end{equation*}
which is only slightly rougher than \eqref{6:apnest}, can be easily derived
from the results of \cite{EHL}.
\vskip0.2cm
As we see it, the techniques developed in the present paper gives a direct
and unified approach to the upper estimates of approximation numbers for
embedding operators of Sobolev spaces on graphs. It can not give any lower
estimates. For the integral operator \eqref{6:ho} our new result consists 
in finding the upper estimate with the best possible constant factor. 

\vskip0.2cm

Some results can be obtained by combination of our both approaches.
For example, the reasonings presented in Subsection 6.2 immediately lead
to the following result.

\begin{prop}\label{6:as} Under the assumptions of \thmref{3:ww}, the
following asymptotic formula for the approximation numbers $a_n$ of the
embedding operator of the space $\tW^{1,p}(\BG,a;o)$ into $\tL^p(\BG,V)$:
\begin{equation*}
\lim_{n\to\infty} n a_n =\a_p\int_\BG w_a(x)V(x)^{1/p}dx,\qquad
\a_p=A_{1,1}([0,1]).
\end{equation*}
\end{prop}

Indeed, for trees this is nothing but a reformulation of Corollary 5.4
from \cite{EHL}. Since the passage to a subspace of finite codimension
does not affect the asymptotic behavior of approximation numbers, the
desired result for general compact graphs follows.

\vskip0.2cm
Lemma 5.9 from \cite{EHL}, which deals with the cases $p=1$ and $p=\infty$,
extends to graphs in the same way.

\subsection{On the case $p\t<1$.} For $\BG=[0,L]$ an analog of
the estimate \eqref{5:11}
follows from \cite{BS2}, Theorem 2.17.  In this analog the assumption
$V\in\tL^1$ is replaced by $V\in\tL^r,\ r=(\t p)^{-1}$, and the term
$\|V\|_1$  in the right-hand side is replaced by $\|V\|_r$. The constant
factor $C(\t,p)$ is of a different nature but still does not depend on
the weight function $V$. 

We do not know whether this result can be extended to arbitrary graphs.
Indeed, it is crucial for our way of reduction the problem to \thmref{2:par}
that the functions $u\in\tL^{\t,p}(\BG)$ are continuous
which is no more true if $p\t<1$.

\bibliographystyle{amsplain}

\begin{thebibliography}{10}


\bibitem{BS1} M.Sh. Birman and and M. Solomyak,
{\it Piecewise polynomial approximations of functions of
classes $W\sb{p}{}\sp{\alpha }$}, Mat. Sb. (N.S.) {\bf 73 (115)}, 1967, 
331--355
(Russian). English transl., Math. USSR-Sb. {\bf 73 (115)} (1967), 295--317.

\bibitem{BS2}  M. Sh. Birman and M. Solomyak, {\it 
Quantitative analysis in Sobolev imbedding theorems and
applications to spectral theory}, 
Tenth Mathem. School, Izd. Inst. Mat.
Akad. Nauk Ukrain  
SSSR, Kiev, 1974, 5--189 (Russian). English transl. in
  Amer. Math. Soc.
Translations, (2), {\bf 114} (1980), 1--132.

\bibitem{ET} D. E. Edmunds and H. Triebel, {\it Function spaces,
entropy numbers and differential operators}, Cambridge tracts in
mathematics 120, Cambridge University Press, Cambridge, 1996.

\bibitem{EH} W. D. Evans and D. J. Harris, {\it Fractals, trees and the
Neumann Laplacian}, Math. Ann. {\bf 296} (1993), 493--527.

\bibitem{EHL} W. D. Evans, D. J. Harris, and J.Lang,
{\it The approximation numbers of Hardy-type operators
on trees}, Proc. London Math. Soc. (3) {\bf 83} (2001) 390--418. 


\bibitem{NS} K. Naimark and M. Solomyak, {\it Eigenvalue estimates for 
the weighted Laplacian on metric
trees}, Proc. London Math. Soc. (3) {\bf 80} (2000), 690--724.
 
\bibitem{P} A. Pinkus, {\it $n$-Widths in Approximation Theory},
Springer-Verlag, Berlin, 1985. 

\bibitem{S} M. Solomyak, {\it On the eigenvalue estimates for the
weighted Laplacian on metric graphs}, in: Nonlinear Problems in
Mathematical Physics and Related Topics I, In Honor of Professor 
O.A. Ladyzhenskaya.

\bibitem{T} H. Triebel, {\it Interpolation Theory. Function Spaces.
Differential Operators}, North-Holland Publishing Co., Amsterdam-New York, 
1978.




\end{thebibliography}

\end{document}